\documentclass[12pt,twoside]{article} 
\usepackage{amsmath} 
\usepackage{amsfonts} 
\usepackage{amsthm} 
\usepackage{amssymb} 
\usepackage{graphicx} 
\usepackage{multicol}

\begin{document} 
\title{{\normalsize{\bf  Fusion of boundary surface-link and ribbonness}}} 
\author{{\footnotesize Akio Kawauchi}\\ 
{\footnotesize{\it Osaka Central Advanced Mathematical Institute, Osaka Metropolitan University}}\\ 
{\footnotesize{\it Sugimoto, Sumiyoshi-ku, Osaka 558-8585, Japan}}\\ 
{\footnotesize{\it kawauchi@omu.ac.jp}}} 
\date\, 

\maketitle 

\vspace{0.25in} 
\baselineskip=10pt 

\newtheorem{Theorem}{Theorem}[section] 
\newtheorem{Conjecture}[Theorem]{Conjecture} 
\newtheorem{Lemma}[Theorem]{Lemma} 
\newtheorem{Sublemma}[Theorem]{Sublemma} 
\newtheorem{Proposition}[Theorem]{Proposition} 
\newtheorem{Corollary}[Theorem]{Corollary} 
\newtheorem{Claim}[Theorem]{Claim} 
\newtheorem{Definition}[Theorem]{Definition} 
\newtheorem{Example}[Theorem]{Example} 

\begin{abstract} A boundary surface-link of at least two components in the 
4-sphere is proved to be a ribbon surface-link if the surface-link obtained from 
it by surgery along a local S-nontrivial fusion 1-handle system is a ribbon surface-link. 
This paper presents a detailed proof of this result, including an appendix,  
although a similar result was previously published with a portion of the proof missing. 
The result is stated using the concept of a \lq\lq local fusion 1-handle\rq\rq, 
which is more precise than the earlier concept of a \lq\lq pairwise fusion 1-handle\rq\rq.

\phantom{x} 


\end{abstract} 

\baselineskip=14pt

\noindent{\bf 1. Introduction} 

\phantom{x} 

This paper presents a detailed version (including an appendix) on the previous result  
\cite[Theorem~1.1]{[1]}, which was published with a portion of the manuscript missing. 
The result is stated using the concept of a \lq\lq local fusion 1-handle\rq\rq, 
which is more precise than the earlier concept of a \lq\lq pairwise fusion 1-handle\rq\rq.
A {\it surface-link of \mbox{$r$} components} is a  closed oriented surface 
$F$ of  $r$ connected components $F_i\, (i =1, 2, \dots, r)$ smoothly embedded 
in the 4-sphere $S^4$. It is called a {\it surface-knot} when $r=1$, 
and a {\it sphere-link} or an 
$S^2$-{\it link} when $F$ consists of 2-spheres. 
A surface-link $F$ of $r(\geq 2)$ components $F_i \, (i = 1, 2, \dots, r)$ 
is a {\it boundary surface-link} if there is a system $V$ of disjoint compact 
connected oriented 3-manifolds $V_i\, (i =1, 2, \dots, r)$ smoothly embedded 
in $S^4$ with $\partial V_i = F_i$, where $V_i$ is  called a {\it Seifert hypersurface} 
of $F_i$. 
The system $V$ is called a {\it disjoint Seifert hypersurface system} 
of the boundary surface-link $F$. 
A {\it 1-handle system} on a surface-link $F$ is a system $h$ of disjoint 1-handles 
$h_j \,(j = 1, 2, \dots, s)$ on $F$ smoothly embedded in $S^4$. 
Let $F(h)$ be the surface-link obtained from $F$ by surgery along a 1-handle system $h$. 
The equivalence on a surface-link is considered setwise, so that the equivalence on the 
surface-link $F(h)$ is independent of a framing of the 1-handle $h$ 
and determined only by the core arc $\alpha$ of $h$ up to 
homotopies of the interior of $\alpha$ in $S^4\setminus F$,  \cite{[2]}. 
For a subsystem $h_1$ of a 1-handle system $h$, let $h_1^c = h \setminus h_1$ 
which induces a 1-handle system on the surface-link $F(h_1)$ as well as 
a 1-handle system on $F$. 
A {\it ribbon surface-link} is the surface-link $O(h)$ obtained from 
a trivial $S^2$-link $O$ by surgery along a 1-handle system $h$, \cite{[3], [4]}.  
A system $h$ of 1-handles $h_j\, (j = 1, 2, \dots, s)$  on 
$F$ is a {\it fusion 1-handle system} if $F(h)$ has just 
$r -s (\geq 1)$ connected components, where each 1-handle  
$h_j$ is  called a {\it fusion 1-handle},  and  a {\it self 1-handle system} 
if  each 1-handle $h_j$ is a self 1-handle, that is, a 1-handle on a connected  
component of $F$.
A special feature of  a boundary surface-link $F$ is the existence 
of a self 1-handle system  $\beta$ on $F$ embedded in a disjoint 
Seifert hypersurface system $V$ for $F$ such that the closed 
complement $V(\beta)=\mbox{cl}(V\setminus \beta)$ is a handlebody system and thus, 
the surface-link $F(\beta)$ is a trivial surface-link in $S^4$, \cite{[2]}. 
For a fusion 1-handle $h$ on a boundary surface-link $F$  
with a disjoint Seifert hypersurface system $V$  in $S^4$, let $F_h$ be 
the surface-link with  two components 
of $F$ to which the 1-handle $h$ attaches, 
and $V_h$ the disjoint Seifert hypersurface system for $F_h$ belonging to $V$.  
Let $F_h^c=F\setminus F_h$ and $V_h^c=V\setminus V_h$. 
The following definition is made for a fusion 1-handle on a boundary surface-link. 

\phantom{x} 

\noindent{\bf Definition.} 
Let $F$  be a boundary surface-link in $S^4$. A fusion 1-handle system $h$ 
on $F$  is {\it local}  if there is a disjoint Seifert hypersurface system $U$ of $F$  
such that every 1-handle belonging to $h$ meets just two components of $U$. 
A local fusion 1-handle $h_1$ on $F$  is  {\it surgically nontrivial},  
or simply {\it S-nontrivial} if $F(\gamma)(h_1)$ is a nontrivial surface-link 
for a self 1-handle system $\gamma$ on $F$ embedded into $U$ and disjoint 
from $h$ such that $U(\gamma)$ is a handlebody system 
for a disjoint Seifert hypersurface system $U$ for $F$ used for the locality of $h_1$. 
Otherwise, $h_1$ is  {\it S-trivial}. A local fusion 1-handle system $h$ on $F$  is 
{\it S-nontrivial} if every 1-handle on $F$  belonging to $h$ is S-nontrivial  
for a disjoint Seifert hypersurface system $U$ of $F$ used for the locality 
of $h$ and  a self 1-handle system $\gamma$ on $F$ embedded into $U$ 
and disjoint from $h$ such that $U(\gamma)$ is a handlebody system. 

\phantom{x} 

The following theorem is the main result. 

\phantom{x}

\noindent{\bf Theorem~1.1.} 
Let $F$ be a boundary surface-link of $r (\geq 2)$ components in $S^4$. If the 
surface-link $F(h)$ for a  fusion 1-handle system $h$ on $F$ consisting of local 
S-nontrivial fusion 1-handles is a ribbon surface-link, then the surface-link $F$  
is a ribbon surface-link.  

\phantom{x}

The proof of Theorem~1.1 is given in Section~4. 
There are many non-ribbon boundary surface-links 
$F$ with two or more components in $S^4$ such that 
the surface-link $F(h)$ for a local {\it S-trivial } fusion 1-handle $h$ 
on $F$ is a ribbon surface-link.  
For example, the  $S^2$-link $L$ consisting of the boundary $S^2$-link 
$K\times \partial I$ of a normal trivial $I$-bundle  $K\times I$ of a non-ribbon 
$S^2$-knot $K$ in $S^4$ and a splitted trivial $S^2$-link $O$ is a  non-ribbon 
boundary $S^2$-link such that the $S^2$-link $L(h)$ for a local S-trivial  fusion 
1-handle $h$ on $L$ obtained by thickening a $I$-fiber  is a trivial $S^2$-link, 
because the 3-manifold $\mbox{cl}(K\times I \setminus h)$ is a 3-ball, where $I=[0,1]$. 
A  surface-link $F$ obtained from this non-ribbon $S^2$-link $L$ by surgery along 
every trivial self 1-handle system disjoint from $h$ is a non-ribbon boundary surface-link  
in $S^4$ such that $h$ is a local S-trivial  fusion 1-handle  on $F$ with $F(h)$  a trivial surface-link, where the non-ribbonness of $F$ is shown from \cite{[5]}. 

{\it Unless otherwise stated, isotopies of \mbox{$S^4$} mean smooth isotopies of 
\mbox{$S^4$} with 4-ball supports, and orientation-preserving diffeomorphisms of \mbox{$S^4$} mean orientation-preserving diffeomorphisms of \mbox{$S^4$} with 
4-ball supports.}

\phantom{x} 

\noindent{\bf 2. Disk basis for trivial surface-link and handlebody system in 
4-sphere} 

\phantom{x} 

For a surface-link $F$ in $S^4$, a {\it loop-pair} on $F$ is a pair of oriented simple loops 
on $F$ of geometric  intersection number $1$, a {\it loop-pair system} on $F$ is a  
system of disjoint loop-pairs on $F$, and a {\it loop basis} of $F$ is a loop-pair system 
on $F$ representing a basis of $H_1(F;Z)$. 
A simple loop $c$ on a surface-link $F$ is {\it spin} if $c$ bounds a possibly 
non-orientable immersed surface $C$ in $S^4$ with 
$C\cap F_1=c$ for the component $F_1$ of $F$ containing $c$ 
such that  the self $Z_2$-intersection number of $C$ in $S^4$ with respect to 
the surface-framing on $F_1$ is $0$, \cite{[6]}.  This is equivalent to saying that 
the $Z_2$-quadratic function $q:H_1(F_1;Z_2)\to Z_2$ 
associated with the surface-knot $F_1$ in $S^4$ sends the 
$Z_2$-homology class of $c$ to the zero.  
Note that the present definition of a spin loop on a surface-link $F$ in $S^4$
 is a generalization of the definition of a spin loop on a surface-knot $F_1$ in 
 $S^4$.  It is a merit of this definition that 
every diffeomorphism $f$ of $S^4$ 
sends a spin loop $c$ on a surface-link $F$ to a spin loop $f(c)$ on  the image 
$f(F)$ of $F$ in $S^4$.
A  loop-pair, a loop-pair system, or a loop basis on a surface-link $F$ is said to be 
{\it spin} if all the loops involved in each case are spin on $F$. 

A meridian loop (namely, the boundary loop of a meridian disk) 
of a handlebody  smoothly embedded in $S^4$ is spin.  
For a loop-pair $(k,\ell)$ on $F$ with $k$ spin on $F$, 
one of the loop-pairs $(k,\ell)$ and $(k,\ell')$ with $\ell'$ 
a Dehn twist of $\ell$ along $k$  on $F$ is a spin loop-pair on $F$.  
Every surface-knot $F_1$ admits a spin loop basis, \cite{[6]}. 
Thus, in the present definition,  {\it every surface-link \mbox{$F$} in \mbox{$S^4$} admits 
a spin loop basis}. 
Every (disjoint) handlebody system $H$ smoothly embedded in $S^4$ is a  split sum 
of  handlebodies of $H$ in $S^4$, meaning that there is a 
disjoint 4-ball system in $S^4$ such that  each component of 
$H$ is contained in a 4-ball belonging to this system, and each 4-ball 
in this system contains exact one component of $H$.
This can be seen by considering $H$ as a spine graph in $S^4$. 
Every trivial surface-link $F$  in $S^4$ is a split sum of 
trivial surface-knots in $S^4$ meaning that there is a disjoint 4-ball 
system in $S^4$ such that each component of $F$ is contained in a 4-ball 
belonging to this system, and each 4-ball in this system contains exact 
one component of $F$. This is because $F$ bounds a handlebody system $H$ in $S^4$. 
Every  handlebody system $H$ smoothly embedded in $S^4$ is  isotopic 
to a standard handlebody system in the standard 3-sphere $S^3$ 
with a smooth isotopy of  $S^4$. Thus, for every handlebody system $H$, 
there is a disjoint disk pair system $(D, E)$ smoothly embedded in $S^4$  with  
$(\partial D, \partial E)=(k,\ell)$ a spin loop basis on the trivial surface-link 
$\partial H=F$ such that $D$ is a meridian disk system of $H$ and $E$ is 
disjoint from the interior of $H$. 
This disk pair system $(D, E)$ is called a {\it disk basis} for the handlebody system 
$H$ in $S^4$. More generally, a disjoint disk pair system 
$(D, E)$ embedded in $S^4$ is called a {\it disk basis} for a trivial surface-link 
$F$ in $S^4$ if $(\partial D, \partial E)=(k,\ell)$  is a spin loop basis on $F$ and 
the interiors of $D$ and $E$ are mutually disjoint and disjoint from of $F$. 
The following {\it uniqueness lemma} on a disk basis for a trivial surface-link is shown.

\phantom{x}

\noindent{\bf Lemma~2.1.}  Let $F$ be a trivial surface-link in $S^4$. 
For every loop basis $(k, \ell)$ on $F$, there is a disk basis $(D,E)$ for $F$ in $S^4$ with 
$(\partial D, \partial E) = (k, \ell)$. For every disk basis $(D,E)$ for $F$, there is a
 handlebody system $H$ bounded by $F$ in $S^4$  such that $(D,E)$ is a disk basis for
  $H$. For any two disk bases $(D,E)$ and $(D', E')$ for $F$, there is 
  an orientation-preserving diffeomorphism $f$ of $S^4$ keeping $F$ 
  setwise and sending $(D,E)$ to 
  $(D', E')$.  
If $(D, E)$ is a disk basis for both handlebody systems $H$ and $H'$ bounded by 
the same $F$ in $S^4$, then $H'$ is moved into $H$ by a smooth isotopy of $S^4$
 keeping $(D,E)$ fixed. 

\phantom{x}

\noindent{\bf Proof.}  
Every trivial surface-link $F$ in $S^4$ admits a disk basis $(D, E)$ 
for a handlebody system $H$ bounded by $F$ since $F$ is moved into a split sum of 
standard surfaces in the standard 3-sphere $S^3$ by a smooth isotopy of $S^4$. 
For any spin loop bases  $(k,\ell)$ and $(k',\ell')$ of a trivial surface-link 
$F$ in $S^4$, there is an orientation-preserving diffeomorphism 
$f$ of $S^4$ with $f'(F)=F$ such that $f'(k, \ell)=(k', \ell')$,  \cite{[7]}.  
Hence for every loop basis $(k,\ell)$ on $F$, there is a disk basis $(D, E)$ for $F$ 
in $S^4$ with $(\partial D, \partial E) = (k, \ell)$, which is also a disk basis for a handlebody system $H$ bounded by $F$ in $S^4$.
For any two disk bases $(D, E)$ and $(D', E')$ for  $F$  with spin loop bases 
$(k, \ell)$ 
and $(k', \ell')$ bounded, respectively, there is an orientation-preserving 
diffeomorphism $f'$ of $S^4$ with $f'(F)=F$ 
 such that the disk bases $(f'(D), f'(E))$ and  $(D', E')$ for $F$ in $S^4$ have 
 the same spin loop basis $(k', \ell')$. Then there is an 
 orientation-preserving diffeomorphism 
 $f''$ 
of $S^4$ with $f''(F) = F$ such that $f''(f'(D), f'(E)) = (D', E')$, \cite{[8], [9]}. 
The composite diffeomorphism 
$f=f''f'$ of $S^4$ has $f(F)=F$ and $f(D,E)=(D', E')$. 
This also means that for every disk basis $(D, E)$ for $F$, there is a handlebody system 
$H$ bounded by $F$ in $S^4$ such that $(D, E)$ is a disk basis for $H$.

Assume that $(D,E)$ is  a disk basis for both handlebody systems $H$ and 
$H'$ bounded by $F$.  Let $H\cup E\times I$ and $H'\cup E\times I$
 be the 3-ball systems 
obtained from  $H$ and $H'$ by adding  the same 2-handle system $E\times I$ 
with core disk system $E$, respectively. 
Let $D\times I$ be the same 2-handle system in $H$ and $H'$ thickening the disk 
system $D$. Then the 3-ball systems $H\cup E\times I$ and $H'\cup E\times I$ 
are moved into each other by 3-cell moves keeping the 3-ball system 
$D\times I\cup E\times I$ fixed. 
Thus, $H'$ is moved into $H$ by an isotopy of $S^4$ keeping $(D,E)$ fixed. 
This completes the proof of Lemma~2.1. 

\phantom{x}

Let $X$ be a bounded smooth 3-submanifold 
of $S^4$, and  $H$ a handlebody system in $S^4$ with $X\cap H=\emptyset$. 
An {\it admissible arc system} between $X$ and $H$  is 
a disjoint simple arc system $\alpha$ whose interior is disjoint from $X\cup H$ 
such that every arc $\alpha_i\, (i=1,2,\dots,s)$ of $\alpha$
joins a point $p_i$  to a point $q_i$ where the points 
$p_i\,(i=1,2,\dots,s)$ are distinct points in $\partial X$ 
and the points $q_i\,(i=1,2,\dots,s)$ are  distinct points
in $\partial H$  belonging  to $s$ distinct components of $H$. 
The following lemma shows well-definedness of a $\partial$-disk sum for a split sum of 
a bounded 3-manifold and a handlebody system.

\phantom{x}

\noindent{\bf Lemma~2.2.}  For any admissible arc system  
$\alpha$ between a bounded smooth 3-submanifold $X$ 
and a handlebody system $H$ in $S^4$  such that $H$ 
is contained in the interior of a 4-ball 
$A$ with 
$A\cap X=\emptyset$, every orientation-preserving diffeomorphism $f$ of $S^4$  
fixing $X$ is $X$-relatively isotopic to a diffeomorphism $f^*$ of $S^4$ fixing 
$X\cup h(\alpha) \cup H$ for a 1-handle system $h(\alpha)$ with 
core arc system $\alpha$. 

\phantom{x}

\noindent{\bf Proof.} 
By considering the 4-ball $A$ as a small 4-ball with 
$A\cap X=\emptyset$, the diffeomorphism $f$ of $S^4$ is 
$X$-relatively isotopic to a diffeomorphism $f'$ of 
$S^4$ with $f'(A)=A$.  Let  $(D,E)$ be a disk basis  for the handlebody system $H$ 
in the interior of the 4-ball $A$ by Lemma~2.1. Then there is an isotopic deformation 
of the 3-ball system $f(H\cup E\times I)$ into the 3-ball system 
$H\cup E\times I$ within the 4-ball $A$  keeping $\partial A$ fixed. 
Thus, the diffeomorphism $f'$ of $S^4$ is $X$-relatively isotopic to a diffeomorphism 
$f''$ of $S^4$ sending the 3-ball system $H\cup  E\times I$ to itself identically. 
By considering the 3-ball system $H\cup  E\times I$ as a finite set regarding 
each 3-ball as a point, the admissible arc system $f''(\alpha)$ is 
$X\cup H$-relatively isotopic to 
$\alpha$ in $S^4$. Hence $f$ is $X$-relatively isotopic to a diffeomorphism $f^*$  of 
$S^4$ fixing $X\cup h(\alpha) \cup H$. This completes the proof of Lemma 2.2.

\phantom{x}

In Lemma~2.2, the 3-manifold $X\cup H\cup h(\alpha)$ is denoted by 
$X \#_{\partial} H$ and called the $\partial$-{\it disk sum} of $X$ and $H$ in $S^4$ 
along the 1-handle system $h(\alpha)$ on $\partial X\cup \partial H$, 
which is determined only by the split sum $X\cup  H$ with the attaching data 
of $\alpha$ and does not depend on choices of the interior of $\alpha$. 
The 1-handle system $h(\alpha)$ is referred to as an {\it admissible} 1-handle system.
A {\it connected sum} of an $S^2$-link $K$ and a trivial surface-link $F$ in $S^4$ 
is a surface-link $K \# F$ obtained from a split sum $K\cup F$ by surgery along 
a 1-handle system $h$ on $K\cup F$ which is  an admissible 1-handle 
system between a Seifert hypersurface $X$ of $K$ and a handlebody system 
$H$ bounded by $F$. 
By Lemma~2.2, the equivalence class of the connected sum $K \# F$ is determined 
only by the split sum $K\cup F$ with the attaching data of the admissible core arc 
system $\alpha$ of  $h$ and does not depend on choices of the interior of $\alpha$.  

\phantom{x}

\noindent{\bf 3. Local fusion 1-handle on  boundary surface-link}  

\phantom{x} 

Let $F$ be a boundary surface-link, and $h$ a local fusion 1-handle on $F$. 
Let $F_h$ be the surface-link with two components of $F$ to 
which the 1-handle $h$ attaches, and $F_h^c = F\setminus F_h$. 
Let $V$ be a disjoint Seifert hypersurface system of $F$, and  $V_h$ 
the disjoint Seifert hypersurface system for $F_h$ belonging to V. Let 
$V_h^c = V\setminus V_h$. Since $h$ is local, take $V$ with 
$h \cap  V_h^c  =\emptyset$. Let $B$ be a 3-ball pair taken from each component 
of the two-component $V_h$ and disjoint from $h$ such that the intersection 
$dL = B \cap  V_h^B$ is a disk pair in the trivial $S^2$-link $L = \partial B$ 
and a proper disk pair in $V_h$, where 
$V_h^B = \mbox{cl}(V_h\setminus B)$.  Let $V^B = V_h^B \cup V_h^c$ 
which is a Seifert hypersurface for the surface-link $F^B = F_h^B \cup  F_h^c$ in 
$S^4$, where 
$F^B = \partial V^B$, $F_h^B = \partial V_h^B$ and $F_h^c = \partial V_h^c$.  
Then $V$ is the $\partial$-disk sum $B \#_{\partial}  V_h^B$ along an admissible 
1-handle pair $hdL$ obtained by thickening and stretching the disk pair $dL$. 
Let $A$ be a 4-ball containing $B$ in the interior and not meeting $V_h^B$ 
whose boundary 3-sphere 
$\partial A$ meets a transverse disk pair of the 1-handle pair  $hdL$. 
The attaching disk pair of the fusion 1-handle $h$ to $F_h$ is slidden into the interior of
 the complementary disk pair $dL^c = \mbox{cl}(L\setminus dL)$. Let $hL$ be a fusion 
 1-handle on $L$ obtained from $h$ by sliding the intersection transverse disk system 
 $h \cap V_h^B$ of $h$ into the interior of B through $V_h^B$, where a thinning
  deformation of $h$ is used. Further, by discarding the $\partial$-disk summand 
  $V_h^B$, the fusion 1-handle $hL$ on $L$ is deformed into the interior of $A$. 
  The following lemma explains some properties on a local S-trivial or S-nontrivial 
  fusion 1-handle on a boundary surface-link. 

\phantom{x} 

\noindent{\bf Lemma~3.1.} Let $h$ be a local fusion 1-handle on a boundary
surface-link $F$ defined on a disjoint Seifert hypersurface system $V$, and 
$h'$ another local fusion 1-handle on $F$ defined on the same $V$ and disjoint from 
$h$ unless $F(h)$ is a surface-knot. 
Let $\beta$ be a self 1-handle system  on $F$ embedded into $V$  and disjoint
from  $h, h'$ such that $V(\beta)$ is a handlebody system. 
The following properties (1)-(4) are obtained.

\phantom{x} 

\noindent{(1)} The local fusion 1-handle $h$ on $F$ is $F(\beta)\cup B$-relatively
 isotopic to the fusion 1-handle $hL$ on $L$ in the interior of the 4-ball $A$ which is
 uniquely determined by $h$ on $F$ up to $F$-relative isotopies of $S^4$, 
and independent of choices of $V$ and $\beta$. 

\phantom{x}

\noindent{(2)} If $h$ is a local S-trivial fusion 1-handle on $F$, then $h$ is 
$F(beta)$-relatively isotopic to a fusion 1-handle $h_0$ on $F$ with 
$h_0\cap V=h_0\cap F$. 

\phantom{x} 

\noindent{(3)} If $h$  is a local S-nontrivial fusion 1-handle on $F$, then  
the surface-link $F(\beta)(h)$ is equivalent to the connected sum 
$L(hL) \# F^B(\beta)$ along the admissible 1-handle pair $hdL$ with 
$L(hL)$ a nontrivial ribbon $S^2$-knot in $A$ and $F^B(\beta)$ 
a trivial surface-link in the 4-ball $A^c = \mbox{cl}(S^4\setminus A)$. 

\phantom{x} 

\noindent{(4)} Unless $F(h)$ is a surface-knot, the surface-link $F(h)$ is always a boundary surface-link and $h'$ is a local fusion 1-handle on $F(h)$. 
Further, $h'$ is S-nontrivial on $F$ if and only if $h'$ is S-nontrivial on $F(h)$.

\phantom{x} 

\noindent{\bf Proof.} 
To prove (1), let $U$ be the disjoint Seifert hypersurface system of $F$ with 
$h\cap  U_h^c = \emptyset$, and $\gamma$  
a self 1-handle system on $F$ embedded into $U$ 
and disjoint from $h$ with $U(\gamma)$ 
a handlebody system, used for the definition of a local S-trivial or 
S-nontrivial fusion 1-handle $h$. Assume that a boundary 
collar of $F_h$ in $U_h$ matches a boundary collar of $F_h$ in $V_h$. 
The $\partial$-disk sums $B \#_{\partial} U_h^B$  and 
$B \#_{\partial} V_h^B$ of $U_h$ 
and $V_h$ are respectively taken along the same admissible 1-handle pair $hdL$. 
The same 4-ball $A$ containing the 3-ball pair $B$ 
but not meeting $U_h^B$ and $V_h^B$ is also taken. 
Let $\alpha$ be the core arc of $h$. Let $\beta_h$ be 
the subsystem of $\beta$ belonging to $V_h$, and 
$\beta _h^c  = \beta  \setminus  \beta_h$.  
By a thinning deformation of handlebody pair $V_h^B(\beta_h)$ 
keeping the meeting data with 
$\alpha$ and keeping $B$ fixed, the handlebody pair $V_h^B(\beta_h)$ meets 
transversely the handlebody pair $U_h^B(\gamma_h)$ with a meridian disk system 
$\mu$ of $V_h^B(\beta_h)$ in the interior of $U_h^B(\gamma_h)$. 
Let $\alpha L^V$ be a simple arc spanning $L$ obtained from $\alpha$ 
by sliding every point of 
the intersection point set $\alpha \cap V_h^B(\beta_h)$ 
into the interior of $B$ through $V_h^B (\beta_h)$, 
and $\alpha^+L$ a simple arc spanning $L$ obtained from 
$\alpha L$ by sliding every point of the intersection point set $\alpha L\cap U_h^B$  
into the interior of $B$ through 
$U_h^B (\gamma_h)$ avoiding the meridian disk system $\mu$ of 
$V_h^B (\beta_h)$ in $U_h^B (\gamma_h)$. Let $\alpha  L^U$ be 
a simple arc spanning $L$ obtained from $\alpha$ by sliding every point 
of the intersection point set 
$\alpha \cap U_h^B (\gamma_h)$ into the interior of $B$ through $U_h^B (\gamma_h)$. 
By these deformations, the arcs $\alpha^+L$ and $\alpha  L^V$ are  
$F(\beta)\cup B$-relatively isotopic to $\alpha$ since they are cell-moves of arcs, 
so that $\alpha^+L$ is $F(\beta )$-relatively isotopic to $\alpha L^V$ in $S^4$. 
By considering the handlebody system $V^B(\beta)$ as a spine graph system, 
the arcs $\alpha^+L$ and $\alpha L^V$ are $V(\beta)$-relatively (and hence 
$F(\beta)\cup B$-relatively) isotopic to $L$-relatively isotopic arcs in the interior 
of $A$, for $\alpha^+L$ and $\alpha L^V$ are disjoint from $V^B(\beta )$ and 
\lq\lq   Regarding $V^B(\beta)$ as a spine graph\rq\rq   is equivalent to 
\lq\lq Discarding $V^B(\beta)$\rq\rq  in this situation. 
Similarly, the arcs $\alpha^+L$ and $\alpha L^U$ are 
$F(\gamma )\cup B$-relatively isotopic to 
$L$-relatively isotopic arcs in the interior of $A$. 
This means that the 1-handle $hL$ with core arc 
$\alpha  L^V$ on $L$ in $A$ is uniquely determined 
up to $L$-relative isotopies in $A$ by the local fusion 1-handle $h$ 
on the boundary surface-link $F$, independently of choices of 
disjoint Seifert hypersurfaces $V$ and self 1-handle systems 
$\beta$ disjoint from $h$ with $V(\beta)$ handlebody systems. 
This completes the proof of (1). 

To prove (2),  assume that $F(\beta)(h)$ is a trivial surface-link. The presentation 
$P=<x_1,x_2|\, x_2^{-1}w x_1 w^{-1}>$  of the fundamental group 
$\pi_1(A\setminus L(hL),p_0)$ of the ribbon $S^2$-knot $L(hL)$ in $A$ is a presentation of an infinite cyclic group and thus there is a generator-preserving isomorphism from 
$P$ to the presentation $<x_1,x_2|\, x_2^{-1}x_1>$,  which implies that a conjugation of the word $x_2^{-1}wx_1 w^{-1}$  is reduced to a word $(x_2^{-1} x_1)^{\pm1}$ in the free group $<x_1,x_2>$ by the conjugacy theorem for the one relator presentation, \cite{[10]}. 
This implies that the arc $\alpha L^V$ is $L$-relatively isotopic in $A$ to a simple arc 
$\alpha_0$ in $A$ not meeting $B$ except for the endpoints, \cite{[11]}. This means that  the 1-handle $h$ is $F(\beta)$-relatively isotpic to 
a 1-handle $h_0$ with the core arc $\alpha_0$ not meeting $V$ except for the attaching part.  This completes the proof of (2).  

To prove (3), assume that the surface-link $F(\beta)(h)$ is a nontrivial surface-link in $S^4$, which means that the presentation$P=<x_1,x_2|\, x_2^{-1}w x_1 w^{-1}>$ of the fundamental group $\pi_1(A\setminus L(hL),p_0)$ of the ribbon $S^2$-knot $L(hL)$ 
in $A$ is a non-abelian presentation by (2), where the generators $x_1, x_2$ are represented by a meridian system of $L$ in $A$. The surface-link $F(\beta)(h)$  is the connected sum $L(hL) \# F^B(\beta)$ of the nontrivial ribbon 
$S^2$-knot $L(hL)$ in $A$ with fundamental group presentation $P$ and the trivial surface-link $F^B(\beta)$ in $A^c$ along the admissible 1-handle pair $hdL$, 
completing the proof of (3).  

To prove (4), assume that $F(h)$ is not a surface-knot. The core arc $\alpha$ of $h$ is deformed to transversely meet the interior 
of $V_h$  with intersection number $0$ by rotating near an endpoint of $\alpha$ 
along a normal circle around $F_h$ in $S^4$. 
Then a Seifert hypersurface $V^*_h$ for $F_h$ with 
$\alpha\cap V^*_h = \alpha\cap F_h$ is obtained from $V_h$ by surgery 
along a 1-handle system whose core arc system belongs to $\alpha$. 
Then $V'= (V_h^*\cup h) \cup V_h^c$ is considered as a Seifert hypersurface 
system for $F(h)$ after a thinning deformation of $h$. Thus, $F(h)$ is a boundary 
surface-link (unless it is a surface-knot). By construction, the 1-handle $h'$ is a local 
fusion 1-handle on $F(h)$ on $V'$ since $h'$ is disjoint from $h$. Further, take a 
$\partial$-disk sum $B'\#_{\partial} {V_{h'}}^{B'}$ of $V$ for a 3-ball pair $B'$ 
taken from each component of the two-component $V_{h'}$ and disjoint from $B$ 
to take a trivial $S^2$-link $L' = \partial B'$ in $S^4$ such that 
$dL'= B'\cap V_{h'}^{B'}$ is a disk pair in $L'$ and a proper disk pair in $V_{h'}$. 
The fusion 1-handle $h'L'$ on $L'$ obtained from $h'$ by using $V$ coincides with the fusion 1-handle $h''L'$ on $L'$ obtained from $h'$ by using $V'$. By (1) and (3), this means that $h'$ is S-nontrivial on $F$ if and only if $h'$ is S-nontrivial on $F(h)$, completing the proof of (4). 
This completes the proof of Lemma 3.1.

\phantom{x} 

In Lemma~3.1 (3), the fundamental group 
$\pi_1(A\setminus L(L^h), p_0)$  of the nontrivial ribbon 
$S^2$-knot $L(hL)$ has a non-abelian presentation 
$P=<x_1,x_2|\, x_2^{-1}w x_1 w^{-1}>$ with generators 
$x_1,x_2$ represented by a meridian system of $L$ in $A$ and with a nontrivial word 
$w$ given by the core arc $\alpha L$
of the 1-handle $hL$ in the free group $\pi_1(A\setminus L, p_0)=<x_1,x_2>$.  
This calculation is classical, \cite{[12]}. 

\phantom{x}

\noindent{\bf 4.  SUPH system for ribbon surface-link and proof of Theorem~1.1}  

\phantom{x} 

A {\it handled sphere system} in $S^4$ is a pair $(O, h)$ embedded in $S^4$ 
such that $O$ is a trivial $S^2$-link, called a {\it based sphere system}, and 
$h$ is a 1-handle system on $O$. 
A ribbon surface-link $F$ is a surface-link constructed from a handled sphere system 
$(O, h)$  by surgery along $h$. 
This pair is  uniquely constructed from a {\it chorded sphere system} 
$(O,\alpha)$ in $S^4$ with a core arc system $\alpha$ of $h$ 
isotopically 
deformed into the standard 3-sphere $S^3$ of $S^4$, 
called a {\it chord system}, or from a {\it chord graph} $(o,\alpha )$ 
in $S^3$ with a trivial link $o$, called a {\it based loop system}, \cite{[3]}. 
The chord graph $(o,\alpha)$ in $S^3$ is considered by a chord diagram 
$C = C(o,\alpha)$ in $S^2$.  
A {\it SUPH system}  for a surface-link $F$ in $S^4$ 
is a compact multi-punctured handlebody system $W$ 
smoothly embedded in $S^4$ whose boundary 
$\partial W$ is given by $\partial W=F \cup O$ for a trivial 
$S^2$-link $O$ in $S^4$. The trivial $S^2$-link $O$ 
is called the {\it based sphere system} 
for the SUPH system $W$.
A SUPH system $W$ for a ribbon surface-link $F$ is 
constructed from  a handled sphere system $(O, h)$  by  taking 
$W= O \times [0,1] \cup h$ for a collar $O\times [0,1]$ of $O$ 
with $O\times \{0\} = O$  in $S^4$, 
where $O\times \{1\} = \partial W \setminus F$ is the based sphere system. Conversely,  
if there is a SUPH system $W$ with based sphere system $O$ 
for a surface-link $F$, then $F$ is a ribbon surface-link constructed from 
a handled sphere system $(O, h)$ with some 1-handle system $h$ in $W$ 
and  every SUPH system $W$ with a fixed based sphere system $O$ 
for a ribbon surface-link $F$ is unique up to boundary collars of $F$ in $W$, \cite{[5]}.
The proof of Theorem~1.1 is done as follows.

\phantom{x} 

\noindent{\bf 4.1: Proof of Theorem~1.1.} 
Let $F$ be a boundary surface-link, and $h$ a local S-nontrivial fusion 
1-handle on $F$ spanning the components $F_i\,(i=1,2)$ of $F_h$. 
By Lemma~3.1 (4),  if the surface-link $F(h)$ is a ribbon surface-link then $F$ 
is a ribbon surface-link, then the proof of Theorem~1.1 
is completed by inductive argument. In fact,  $F(h)$ is a boundary surface-link 
(unless it is a surface-knot) and the fusion 1-handle 
$h''$ on $F(h)$ induced from a next local S-nontrivial fusion 1-handle $h'$ on $F$ is 
a local S-nontrivial fusion 1-handle on $F(h)$. 
Assume that $h\cap V_h^c=\emptyset$ for a disjoint Seifert hypersurface $V$ of $F$. 
Let $B$ is a 3-ball pair taken from each component of 
the two-component $V_h$ and disjoint from $h$ 
such that  the intersection $dL=B\cap V_h^B$ is a 
disk pair in the trivial $S^2$-link $L = \partial B$ 
and a proper disk pair in $V_h$, where 
$V_h^B=\mbox{cl}(V_h\setminus B)$. 
By Lemma 2.2,  $V$ is taken as the $\partial$-disk sum 
$B \#_{\partial}  V^B$ along an admissible 1-handle pair $hdL$ 
obtained by thickening and stretching the disk pair $dL$. 
The 3-ball pair $B$ is taken in the interior of a 4-ball $A$ 
not meeting $V^B$ whose boundary $\partial A$ meets $hdL$ 
with a parallel copy of $dL$. 
Let $\beta$ be a self 1-handle system on $F$ embedded in $V^B$ 
and disjoint from $h$ such that $V^B(\beta)$ is a handlebody system, 
whose boundary is a trivial surface-link $F^B(\beta)$. 
The following assertion is obtained by Lemma 3.1 (1) and (3).

\phantom{x} 

\noindent{\bf (4.1.1)} 
The $\partial$-disk sum 
$W'= W(L;hL) \#_{\partial} V^B(\beta)$ along an admissible 1-handle pair $hdL$ 
 obtained by thickening and stretching the disk pair $dL$ is a SUPH system for the surface-link $F(\beta)(h)$ such that $W(L;hL)$ is 
the union SUPH system $B^{(0)} \cup hL$ for the nontrivial ribbon $S^2$-knot 
$L(hL)$ in the interior of $A$ and $V^B(\beta)$
is the handlebody system  bounded by $F^B(\beta)$ in the interior of the 4-ball 
$A^c = \mbox{cl}(S^4\setminus A)$ where  $B^{(0)}$ is a boundary collar  of $L$ in 
$B$ and $hL$ is a 1-handle on $L$ constructed from $h$. 

\phantom{x}  

By assumption, the surface-link $F(h)$ is a ribbon surface-link. Let $W(F(h))$ be 
a SUPH system for $F(h)$. 
Replace $W(F(h))$ with a compact multi-punctured manifold 
$W(F(h))^{(0)}$ of $W(F(h))$ to avoid the intersection with $\beta$  
except for the attaching part of $\beta$ 
and to avoid having an unpunctured handlebody component. 
Then the union $W = W(F(h))^{(0)} \cup \beta$ is a SUPH system 
without unpunctured handlebody components for $F(\beta)(h)$. 
Then the following assertion is obtained.

\phantom{x} 

\noindent{\bf (4.1.2)} 
There is an orientation-preserving diffeomorphism $f$ of $S^4$  
sending the SUPH system $W'$ for $F(\beta)(h)$ to a SUPH system $W_{**}^+$ for 
$F(\beta)(h)$ obtained from some multi-fusion SUPH system 
$W_{**}$ of some multi-fission SUPH system $W_*$ 
of $W$ by adding some 3-ball system $B^W$ along some 
based sphere system of $W_{**}$.

\phantom{x} 

\noindent{\bf Proof.} Let ${W'}^{(0)}$ be a multi-punctured SUPH system of $W'$ 
to avoid unpunctured handlebody components. By Theorem A.1 of Appendix, there is an orientation-preserving diffeomorphism $f$ of $S^4$  sending a SUPH system 
${W'}^{(0)}$ for $F(\beta)(h)$ to some $W_{**}$, which sends the SUPH system 
$W'$ for $F(\beta)(h)$ to some $W_{**}^+$. This completes the proof of (4.1.2).  

 \phantom{x}
 
The S-nontrivial fusion 1-handle $hL$ meets the 3-ball pair $B$ 
with transverse disks by Lemma~3.1 (3). Let $D_h$ be such a transverse disk of $h$, 
and $\partial D_h = k_h$.    
The proper disk $D_h^f = f(D_h)$  in $W_{**}^+$ is assumed to be 
$D_h^f \cap B^W = \emptyset$, so that $D_h^f$ is in $W_{**}$ and hence in $W_*$. 
It is considered that the diffeomorphism $f$ sends the $\partial$-disk sum 
$B^{(0)} \#_{\partial} V^B(\beta)$ along $hdL$ which is a SUPH system for $F(\beta)$  
to the SUPH system $R = W_{**}^+(D_h^f\times I)$ for $F(\beta)$, 
and  that  $h$ and $h^f = f(h)$ are $F(\beta)$-relatively isotopic 1-handles on 
$F(\beta)$ in $S^4$ with the same attaching disk pair in $F(\beta)$, by Lemma 3.1 (1). 
Let $D_{\beta}$ be a transverse disk 
system of the 1-handle system $\beta$ in $W$ with just one disk for every 1-handle 
of $\beta$, and $\partial D_{\beta} = k_{\beta}$. 
The disk $D_h^f$  has $\partial D_h^f = k_h$ disjoint from $k_{\beta}$ and 
may meet the interior of the disk system $D_{\beta}$ in $W_* $ transversely.  
A proper disk system $D_h^f*$ in $W_* $ with $D_h^f* \cap D_{\beta} = \emptyset$ 
and $\partial D_h^f* = k_h$ is constructed from $D_h^f$ in $W_* $ by cut and paste operations on the transverse intersection loop system $D_h^f\cap D_{\beta}$ in the disk $D_{\beta}$. By these constructions, there is a finite sequence 
\[D_h^f* = D_0 \to D_1 \to \dots \to D_n  = D_h^f\] 
of proper disks $D_i$ with $\partial D_i = k_{\beta}$ in $W_*$ 
such that  $D_h^f* \cap D_{\beta} = \emptyset$ and $D_{i+1} = \mbox{cl}(D_i\setminus \delta_i)\cup \delta'_i$  ($i = 0, 1, \dots n-1$), where $\delta_i$ and $\delta'_i$ are disks with $\partial \delta_i = \partial \delta'_i$ 
such that the union $K_i = \delta_i \cup \delta'_i$ is a ribbon $S^2$-knot in 
the interior of $W_*$ with a SUPH system $U_i$ in $W_*$ for every $i$.  
Let $W_*\times [0, 1]$ be a collar of $W_*$ in $S^4$  with $W_*\times 0 = W_*$. 
Let $0 = t_0 < t_1 < · · · < t_n = 1$ be a division of the interval $[0, 1]$.
In $W_*\times [0, 1]$, let $D_i[t_j ] = k_h\times [0, t_j ] \cup D_i\times t_j$ and
$U_i[t_j ] = U_i\times t_j$.
Let  $D_0[0]\times I = D_0 \times I$ be a 2-handle on $F(\beta)(h)$ 
thickening $D_0$  in $W_*$, and 
$D_i[t_j ] \times  I = k_h\times  [0, t_j] \times  I \cup D_i\times  t_j \times  I$ 
 be a 2-handle on $F(\beta)(h)$ thickening $D_i$ in $W_* \times [0, 1]$. 
 Let $F_i = F(\beta)(h)(D_{\beta} \times  I)(D_i[t_i]\times  I)$ $(i = 0, 1, \dots, n)$ 
 and $F'_{i+1} = F(\beta)(h)(D_{\beta} \times I)(D_i[t_{i+1}]\times I)$  
 $(i = 0, 1, \dots, n-1)$. 
The surface-link $F_0$ is a ribbon surface-link since the compact 3-manifold 
$W_* ^-$ obtained from $W_* $ by cutting along the disk system 
$D_h^f* \cup D_{\beta}$ is a SUPH system for $F_0$. 
Then $F'_1$ is a ribbon surface-link since it is equivalent to $F_0$. 
The surface-link $F_1$ is a ribbon surface-link by Theorem~A.2 of Appendix, 
because the based sphere system of $U_i$ belongs to the based sphere system of 
$W_*$ and $F_1$ is obtained from $F'_1$ as twice SUPH moves on $U_1[t_1]$. 
Next, $F'_2$ is a ribbon surface-link since it is equivalent to $F_1$. 
The surface-link $F_2$ is a ribbon surface-link since it is obtained from $F'_2$ as twice
 SUPH moves on $U_2[t_2]$. By continuing this procedure, the surface-link 
$F_n = F(\beta)(h)(D_{\beta}\times I)(D_n[t_n]\times I) 
= F(\beta)(h)(D_{\beta}\times I)(D_h^f[1]\times I)$ 
is a ribbon surface-link. Since the 1-handles $h$ 
and $h^f$ on $F(\beta)$ have the same attaching disk pair in $F(\beta)$, 
the surface-link $F_n$ is equivalent to 
$F(\beta)(h)(D_{\beta}\times I)(D_h \times I) =F$ by considering 
the 1-handle 
$h^f[1]= (h^f \cap F) \times  [0,1] \cup h^f\times 1$ 
in $W_* \times [0, 1]$ instead of $h^f$. 
Thus, the surface-link $F$ is a ribbon surface-link. 
This completes the proof of Theorem 1.1. 
 
\phantom{x}

\noindent{\bf  Appendix: Moves on the SUPH systems of a ribbon surface-link}

\phantom{x} 

There are three kinds of moves $M_0, M_1, M_2$ on a chord diagram $C = C(o,\alpha )$ 
for equivalence of  a ribbon surface-link $F$, which are called the 
{\it Reidemeister move}, the {\it fusion-fission move}, and  the {\it chord move}, 
 respectively, \cite{[3],[13]}. 
The fusion-fission move $M_1$ consists of the fusion move $M_1(fusion)$ 
decreasing the number of based loops and the fission move $M_1(fission)$ 
increasing the number 
of based loops. 
For a SUPH system $W'$ obtained from a SUPH system $W$ by applying 
the Reidemeister move $M_0$ or the chord move $M_2$,  there is a  
diffeomorphism of $S^4$ 
isotopic to  to the identity and sending  $W$ into $W'$.  Namely, $W'$ 
is an isotopic deformation of $W$.
If the fusion move  $M_1(fusion)$ or its iteration produces a SUPH system $W'$ 
 from a SUPH system $W$, then $W'$ is obtained from $W$ by surgery of the based 
 sphere system $O$ of $W$ along a fusion 1-handle or a fusion 1-handle system 
  embedded in $W$, which is called a {\it fusion SUPH system} or a {\it multi-fusion 
 SUPH system} of  $W$, respectively. The fusion move $M_1(fusion)$ and the 
 fission move $M_1(fission)$ are geometrically dual concepts on the SUPH systems 
 for a ribbon surface-link.  
Thus, if the fission move $M_1(fission)$ or its iteration produces 
a SUPH system $W'$ from a SUPH system $W$, then $W'$ is 
obtained from $W$ by attaching a 2-handle or a 2-handle system on 
the based sphere system $O$ of $W$, called a {\it fission SUPH system} or 
a {\it multi-fission SUPH system} of $W$, respectively, where note that 
if a 2-handle $h^2$ on a component $O_1$ of $O$ is used for 
a fission SUPH system 
$W'$ of the SUPH system $W$, then there is a 3-ball $B_1$ with 
$\partial B_1 = O_1$ in $S^4$ such that the intersection 
$d_1 = W \cap B_1$ is a  2-disk system in the interior of 
$B_1$ (coming from the 1-handle system of a handled sphere system 
of $W$) and the 2-handle $h^2$ on $O_1$ is embedded in $B_1$ with 
$h_2 \cap d_1 = \emptyset$.  Any two SUPH systems $W$ and $W'$ for a ribbon 
surface-link $F$ excluding handlebodies are moved into each other 
by a finite number of the moves $M_i (i=0,1,2)$, \cite{[13]}. 
Up to up to isotopic deformations, only a finite number of the fusion move $M_1(fusion)$ 
and/or the fission move $M_1(fission)$ remain. The following theorem is shown.

\phantom{x} 

\noindent{\bf Theorem~A.1.} Let $W$ and $W'$ be SUPH systems  of a ribbon 
surface-link $F$ without handlebody components.  
Then a multi-fusion SUPH system $W_{**}$ of a multi-fission SUPH system 
$W_*$ of $W$ up to isotopic deformations is isotopic to the SUPH system $W'$ in 
$S^4$. In other words, a multi-fission SUPH system $W_*$ of  $W$ up to isotopic deformations is isotopic to a multi-fission SUPH system $W'_*$ of $W'$ in $S^4$.  

\phantom{x}

\noindent{\bf Proof.} 
Assume that a SUPH system $W'$ is obtained from a SUPH system $W$ by 
$M_1(fusion)$ and a SUPH system $W''$ is obtained from $W'$ by $M_1(fission)$. Then 
both $W$ and $W''$ are obtained from $W'$ by $M_1(fission)$. After an isotopic 
deformation of $W'$ in $S^4$, assume that the 2-handles $h^2$ and ${h^2}''$ used to 
produce $W$ and $W''$ from $W'$ are embedded to the same 3-ball system $B'$ 
bounded by the basic sphere systems $O'$ of $W'$ although the 1-handle system of 
the handled sphere system of $W'$ may have a different position for $W$ and $W''$. 
The 2-handles 
$h^2$ and ${h^2}''$ are made disjoint by an isotopy of $B$ (while paying attention to the
intersection between $B$ and the 1-handle system of the handled sphere system) even  
if they belong to the same component of $B$. Then $W'\cup h^2$ is isotopically 
deformed to 
$W$ and  $W'\cup {h^2}''$ is isotopically deformed to $W''$. The move from 
$W'\cup h^2$ to $W'\cup h^2\cup {h^2}''$ is $M_1(fission)$ and the move from 
$W'\cup h^2\cup {h^2}''$ to $W'\cup {h^2}''$ is $M_1(fusion)$. 
This means that the ordering of the moves of $M_1(fusion)$ 
and $M_1(fission)$ is arranged up to isotopic deformations. 
Since the fusion move $M_1(fusion)$ and the fission move 
$M_1(fission)$ are geometrically dual concepts, a multi-fission SUPH system 
$W_*$ of $W$ 
up to isotopic deformations is isotopic to a multi-fission SUPH system $W'_*$ 
of $W'$. This completes the proof of Theorem A.1.  

\phantom{x}

Let $F$ be surface-linkin $S^4$, and $U$  a SUPH system for a ribbon $S^2$-knot  $K$ 
 in $S^4$ such that $F\cap U=\emptyset$.
The surface-link $F'=(F\cup K)(\eta)$ for  a fusion 1-handle $\eta$ on $F\cup K$ 
with $U\cap \eta=K\cap \eta$ a disk is called a {\it SUPH move} of the surface-link 
$F$ by a SUPH system  $U$.  Let $O^U$ be the based sphere system of $U$. 
The following proposition is used for the proof of Theorem 1.1. 

\phantom{x} 

\noindent{\bf Theorem~A.2.} The  SUPH-move $F'$ of a ribbon surface-link $F$ with  
based sphere system $O$ by a SUPH system  $U$ is a ribbon surface-link if  the union 
$O\cup O^U$ is a trivial $S^2$-link.

\phantom{x} 
 
\noindent{\bf Proof.} The  ribbon surface-link $F$ is given by a handled sphere system 
$(O,h)$. Let $(O^U,h^U)$ be a handled sphere system for a ribbon $S^2$-knot $K$  
given by $U$, where $h$, $h^U$ and $\eta$ are made mutually disjoint for a 
fusion 1-handle $\eta$ joining $F$ to $K$. 
The  SUPH-move $F'$ of $F$ is obtained as the ribbon surface-link 
$(O\cup O^U)(h\cup h^U\cup \eta)$. 
This completes the proof of Theorem~A.2. 
 
\phantom{x}

\noindent{\bf Acknowledgements}  

The author thanks the reviewers for pointing out some inaccuracies of shortened 
 explanations in the older version of this paper.
This work was partly supported by JSPS KAKENHI Grant Numbers JP21H00978 
and JP26K06456, and MEXT Promotion of Distinctive Joint Research Center 
Program JPMXP0723833165 and Osaka Metropolitan University Strategic 
Research Promotion Project (Development of International Research Hubs). 

\phantom{x}

\end{document}